\documentclass[12pt]{article}

\usepackage{latexsym}
\usepackage{amsmath}
\usepackage{amssymb}
\usepackage{amsfonts}

\begin{document}
\parskip 0.2 cm
    \renewcommand{\theequation}{\arabic{equation}}

\def \o{\over}
\def \s{\sigma}
\def \G{\Gamma}
\def \g{\gamma}
\def\gm{\gamma(dy,du)}
\def\xm{\xi(dy,du)}
\def \a{\alpha}
\def \l{\lambda}
\def \disp{\displaystyle}
\def \up{\uparrow}
\def \U{{\cal U}}
\def \K{{\cal K}}
\def \P{{\cal P}}
\def \D{{\cal D}}
\def \M{{\cal M}}
\def \Mp{{\cal M_+}}
\def \ve{\varepsilon}
\def\reals{I\!\!R}
\def \O{\Omega}
\def \ph{\varphi}
\def \g{\gamma}
\def \t{\tau}
\def\hf{\hfill{$\Box$}}
\def \dn{\downarrow}
\def \up{\uparrow}
\def \d{\delta}
\def \om{\omega}

\newtheorem{Theorem}{Theorem}[section]
\newtheorem{Proposition}[Theorem]{Proposition}
\newtheorem{Remark}[Theorem]{Remark}
\newtheorem{Lemma}[Theorem]{Lemma}
\newtheorem{Corollary}[Theorem]{Corollary}

\begin{center}
{\bf \Large{Lack of Equality between Abel and Ces\`aro Limits in Discrete Optimal Control and the Implied Duality Gap}}
 \bigskip
 \bigskip

 {\bf  Ilya Shvartsman}\\
{\it \small{Department of Mathematics and Computer Science,
Penn State Harrisburg, Middletown, PA 17057, USA}}
\end{center}

\bigskip
{\bf Abstract.} In the recent paper \cite{BGS} it has been shown that if  Ces\`aro and Abel limits for a certain discrete time optimal control problem are not equal, then there is a 
duality gap between a certain infinite-dimensional linear programming problem and its dual. In this paper we construct an example of a problem satisfying the assumptions of \cite{BGS} where  Ces\`aro and Abel limits are different.



\section{Introduction}

Consider the discrete time controlled dynamical system
\begin{equation}\label{A1}
\begin{aligned}
&y(t+1)=f(y(t),u(t)), \; t\in N=\{0,1,\dots\}\, \\
&y(0)=y_0,\\
&y(t)\in Y,\\
&u(t)\in U(y(t)).
\end{aligned}
\end{equation}
Here  $Y$ is a subset of $\reals^m$, $\ U(\cdot):\,Y\leadsto U_0$ is an upper semicontinuous compact-valued mapping to a  compact metric space $U_0$,
$\ f(\cdot, \cdot ):\,\reals^m\times U_0\to \reals^m$ is a continuous function.

On  the trajectories of \eqref{A1}, we consider the following optimal control problems:
\begin{equation}\label{A112-1}
{1\o T} \min_{u(\cdot)}\sum_{t=0}^{T-1} g(y(t),u(t))=:V_T(y_0),
\end{equation}
\begin{equation}\label{A112-2}
(1-\a) \min_{u(\cdot)}\sum_{t=0}^{\infty}\a^t g(y(t),u(t))=:h_{\a}(y_0),
\end{equation}
where $g:\,\reals^m\times U_0\to \reals^m$ is a continuous function and $\a\in (0,1)$ is a discount factor. The minimima are taken with respect to controls satisfying the constraints of \eqref{A1} on a finite horizon $t\in \{0,\dots,T\}$ in \eqref{A112-1} and on infinite horizon in  \eqref{A112-2}. 
The former is a problem of minimizing of the long-run average and the latter is a problem with discounting.

The limits
\begin{equation}\label{A2}
 \lim_{T\rightarrow\infty}V_T(y_0) \hbox{ and } \lim_{\a\up 1}h_{\a}(y_0) \,,
\end{equation}
are often referred to as {\em Ces\`aro} and {\em Abel limits}, respectively. Such limits have been studied in various contexts  since the work of Hardy and Littlewood \cite{HardyLit14}.  They proved  that if $\{a_n\}_{n\ge 1}$ is a bounded sequence of real numbers, the existence of one of the limits 
$$
\lim_{n\to\infty} {1\o n}\sum_{i=1}^{n} a_i \hbox{ and } \lim_{\a\up 1} (1-\a)\sum_{i=1}^{\infty} \a^i a_i
$$ 
implies the existence of the other and their equality. (The equality of  Ces\`aro and Abel limits is often referred to as {\em Tauberian theorem}.) In \cite{Hardy49} a similar result was proved for real-valued functions: if $g$ is countinuous and bounded on the real line, then the existence of one of the limits 
$$
\lim_{T\to\infty} {1\o T}\int_{0}^{T} g(t)\,dt \hbox{ and } \lim_{\l\dn 0}  \l\int_{0}^{\infty} e^{-\l t}g(t)\,dt
$$ 
implies the existence of the other and their equality.

There is an extensive literature devoted to the existence and equality of Ces\`aro and Abel limits in problems of dynamic programming and optimal control in discrete and continuous time, see, e.g., \cite{Arisawa-3,Bardi,BQR-2015,GQ,GruneSIAM98,GruneJDE98,Khlopin,Sorin92,OV-2012,QR-2011}.  In \cite{Sorin92} and \cite{OV-2012} it is shown that Tauberian theorem holds in discrete and continuous time, respectively, if convergence in  Ces\`aro and Abel limits is {\em uniform} with respect to the initial state. There is an example in \cite{Sorin92} in discrete time and discrete unbounded state space,  where both  Ces\`aro and Abel limits exist, but are not equal to each other. Papers \cite{OV-2012} and \cite{QR-2011} exhibit an example in continuous time with a compact state constraint set $Y$ where this equality is lacking. It is worth mentioning that the equality 
\begin{equation}\label{A22}
 \lim_{T\rightarrow\infty}\min_{y_0\in Y}V_T(y_0) =\lim_{\a\up 1}\min_{y_0\in Y}h_{\a}(y_0) \,,
\end{equation}
where minimization taken with respect to the initial condition, is always valid (\cite{GPS-2017}).



In the recent paper \cite{BGS} it has been shown that for problems \eqref{A112-1} and \eqref{A112-2} with a compact set $Y$, the Ces\`aro and Abel limits are bounded above
by the optimal value of a certain infinite-dimensional linear programming (IDLP) problem and
bounded from below by the optimal value of the corresponding dual IDLP
problem (see Proposition 2.3 and Theorem 3.1 in \cite{BGS}). This implies that the Ces\`aro and Abel limits must exist and be equal if the strong duality holds, that is, the optimal values of the IDLP and its dual are equal. 

While in the finite-dimensional case strong duality always holds, it may not be true in linear programming in infinite dimensions (see examples in \cite{AN}). In \cite{Gonzalez98} it is shown that there is no duality gap for a certain class of IDLP problems in stochastic setting. In this paper we construct a discrete-time example of a system satysfying the assumptions of \cite{BGS} where the limits \eqref{A2} exist, but are different. This example is importnat not only for its own sake, but it also demonstrates the presense of a duality gap. Ideas from \cite{Sorin92}, \cite{OV-2012} and \cite{QR-2011} were used in the contruction of this example.



\section{Example}

In this section we demonstrate a step-by-step construction of the example, starting from a motivating contintinous time uncontrolled system through a discrete-time controlled system. To be consistent with the assumptions of  Theorem 4.1 in \cite{BGS}, we need to ensure that the set $Y$ is compact,  the function $V_T(\cdot)$ is continuous for all $T$ and the function $h_{\a}(\cdot)$ is continuous for all $\a\in (0,1)$.

\subsection{Continuous Time Uncontrolled System}

Consider the following uncontrolled system in continuous time
\begin{equation}\label{VV1}
\begin{aligned}
&\dot x=y, \; t>0,\\
&\dot y=0,\\
&(x(0),y(0))=(x_0,y_0).
\end{aligned}
\end{equation}
For $y_0>0$ the trajectory moves along a horizontal line $y=y_0$ with a constant speed $y_0$. If $y_0=0$, there is no motion.
For natural $T$ and $\l>0$ define 
$$
\g^c_T(x_0,y_0)={1\o T}\int_0^T g(x(t))\,dt
$$ 
and
$$
v^c_{\l}(x_0,y_0)=\l\int_0^{\infty} e^{-\l t}g(x(t))\,dt, 
$$
where
$$
g(x)=\begin{cases} 0,&1\le x\le 2,\\
1,& 0\le x<1,\; x>2. \end{cases}
$$
(This $g$ is discontinuous, but is can be smoothed out, see Section 2.3. Also, the dynamics of the system can be adjusted so that it evolves on invariant compact set, see the end of Section 2.4.)

Denote $D:=[0,1)\times [0,1]$.

{\bf Proposition 1.} {\em
We have  $\g^c_T(x_0,y_0)\to 1$ as $T\to \infty$ for all $(x_0,y_0)\in D$ and for all $x_0\in [0,1]$
\begin{equation}\label{VV5}
\min_{y_0\in [0,1]}\g^c_T(x_0,y_0)= {1-x_0\o 2-x_0}\,.
\end{equation}
}

{\bf Proof.} The convergence  $\g^c_T(x_0,y_0)\to 1$ as $T\to \infty$ follows from the fact that any trajectory ultimately spends almost all time in the region where $g=1$.
Formula \eqref{VV5}
follows from the fact that the minimum of $\g^c_T(x_0,\cdot)$ is achieved when the trajectory is such that $x(T)=2$, since in this case the trajectory $x(t)$ spends maximum possible proportion of time in the band $1<x<2$ where $g=0$. This occurs when $y^*_0=(2-x_0)/T$ and $\disp \g^c_T(x_0,y^*_0)={1\o T}{(1-x_0)\o y^*_0}=
{1-x_0\o 2-x_0}$. (Notice that $T\ge 2$ ensures that $y_0^*\in [0,1]$.)

\hf

Let us analyse the function $v^c_{\l}(x_0,y_0)$. 

{\bf Proposition 2.} {\em

{\rm (a)} For all $x_0\in [0,1)$
 \begin{equation}\label{VV2}
v^c_{\l}(x_0,y_0)=\begin{cases}
1-e^{-{\l(1-x_0)/y_0}}+e^{-{\l(2-x_0)/y_0}}\,,\, &y_0\in (0,1]\\
1, &y_0=0.
\end{cases}
\end{equation}
Consequently  $v^c_{\l}(x_0,y_0)\to 1$ as $\l\to 0$ on $D$.

{\rm (b)} For $x_0\in [0,1)$ and $\l<\ln 2$ 
\begin{equation}\label{VV5.5}
\min_{y_0\in [0,1]}v^c_{\l}(x_0,y_0)=1-{(1-x_0)^{1-x_0}\o (2-x_0)^{2-x_0}}\,.
\end{equation}
}

{\bf Proof.} For $y_0\in (0,1]$ and $x_0\in [0,1)$ the values of $t$ when the trajectory reaches $x=1$ and $x=2$ are equal to 
\begin{equation}\label{t_1}
t_1:={1-x_0\o y_0}\quad\hbox{and }t_2:={2-x_0\o y_0},
\end{equation}
respectively.
Therefore,
$$
v^c_{\l}(x_0,y_0)=\l\left(\int_0^{(1-x_0)/y_0} e^{-\l t}\,dt+\int_{(2-x_0)/y_0}^{\infty} e^{-\l t}\,dt\right)=1-e^{-\l(1-x_0)/y_0}+e^{-\l(2-x_0)/y_0},
$$
which is part (a) of the proposition. 

Denoting $s=e^{-\l/y_0}$ we have 
\begin{equation*}
v^c_{\l}(x_0,y_0)=1-s\cdot s^{-x_0}+s^2\cdot s^{-x_0}=1-s^{1-x_0}(1-s):=\ph(s), 
\end{equation*}
hence,
\begin{equation}\label{BB1}
\inf_{y_0\in (0,1]}v^c_{\l}(x_0,y_0)=\inf_{s\in (0,e^{-\l}]}\ph(s). 
\end{equation}
Since $y_0=0$ and $s=0$ are not minimizers of $v^c_{\l}(x_0,\cdot)$ and $\ph(\cdot)$, respectively, \eqref{BB1} implies that
\begin{equation}\label{min_v}
\min_{y_0\in [0,1]} v^c_{\l}(x_0,y_0)=\min_{s\in [0,e^{-\l}]}\ph(s).
\end{equation}
We have $\ph'(s)=s^{-x_0}(x_0-1+(2-x_0)s)$, hence, the minimim of $\ph$ is reached when  
$s=(1-x_0)/(2-x_0)$. (Note that this $s$ lies in the interior of the interval $[0,e^{-\l}]$ for $\l<\ln 2$ for all $x_0\in[0,1)$). Substituting this value of $s$ into the right side of \eqref{min_v} we obtain \eqref{VV5.5}. \hf


\subsection{Discrete Time Uncontrolled System}

Consider a discrete-time counterpart of system \eqref{VV1}:
\begin{equation}\label{VV4}
\begin{aligned}
&x(t+1)=x(t)+y(t), \; t\in N,\\
&y(t+1)=y(t),\\
&(x(0),y(0))=(x_0,y_0).
\end{aligned}
\end{equation}

Denote the corresponding costs by
$$
\g^d_T(x_0,y_0)={1\o T}\sum_{t=0}^{T-1} g(x(t))
$$ 
and
$$
v^d_{\a}(x_0,y_0)=(1-\a)\sum_{t=0}^{\infty} \a^{t}g(x(t)), 
$$
where $\a\in(0,1)$.

{\bf Proposition 3.} {\em 
For $x_0\in[0,1)$ we have
\begin{equation}\label{ZX1}
\inf_{y_0\in [0,1]}\g^d_T(x_0,y_0)= {1-x_0\o 2-x_0}+\eta_1(T),
\end{equation}
where $\eta_1(T)\to 0$ as $T\to \infty$ uniformly with respect to $x_0$.
}

\bigskip

{\bf Proof.} Notice that due to discontinuity of $g$, existence of minimumin in \eqref{ZX1} can't be ascertained.
Due to the specific piecewise-constant structure of $g$, for $y_0\in (0,1]$ and $x_0\in [0,1)$, we have 
\begin{equation}\label{KO1}
|\g^d_T(x_0,y_0)-\g^c_T(x_0,y_0)|\le {2\o T}\,.
\end{equation}
It is also clear that  $\g^d_T(x_0,0)-\g^c_T(x_0,0)=0$ for all $T$ and $x_0\in [0,1)$, hence, \eqref{KO1} holds for all $y_0\in [0,1]$. 

For arbitrary $\d>0$ let $\tilde y$ be such that  $\g^d_T(x_0,\tilde y)\le \inf_{y_0\in [0,1]}\g^d_T(x_0,y_0)+\d$ and
let $y_0^*$ be the minimizer of $\g^c_T(x_0,\cdot)$, as in Proposition 1. We have
$$
\inf_{y_0\in [0,1]}\g^d_T(x_0,y_0)-\min_{y_0\in [0,1]}\g^c_T(x_0,y_0)\le \g^d_T(x_0,y_0^*)-\g^c_T(x_0,y_0^*)\le {2\o T}
$$
and 
$$
\inf_{y_0\in [0,1]}\g^d_T(x_0,y_0)-\min_{y_0\in [0,1]}\g^c_T(x_0,y_0)\ge\g^d_T(x_0,\tilde y)-\d-\g^c_T(x_0,\tilde y)\ge -{2\o T}-\d.
$$
From the fact that $\disp \min_{y_0\in [0,1]}\g^c_T(x_0,y_0)= {1-x_0\o 2-x_0}$ due to Proposition 1 and the inequalities above, we conclude that \eqref{ZX1} is true
with $\eta_1(T)=2/T$.

\hf

{\bf Proposition 4.} {\em 
For all $x_0\in[0,1)$ and $\a\in[1/2,1)$ we have
\begin{equation}\label{VV5.7}
\inf_{y_0\in [0,1]}v^d_{\a}(x_0,y_0)=1-{(1-x_0)^{1-x_0}\o (2-x_0)^{2-x_0}}+\nu_1(\a)\,,
\end{equation}
where $\nu_1(\a)\to 0$ as $\a\up 1$ uniformly with respect to $x_0$.
}

{\bf Proof.} 
For $x_0\in [0,1)$ and $y_0\in(0,1]$ let $\t_1\in N$ be the maximum time when the trajectory of \eqref{VV4} satisfies $x(\t_1)\le 1$, and $\t_2\in N$ be the minimum time when this trajectory satisfies $x(\t_2)\ge 2$. We have $\t_1=[t_1]$ and
\begin{equation}\label{VV11}
\begin{aligned}
\t_2=\begin{cases} t_2 &\hbox{ if }t_2 \hbox{ is an integer}\\
[t_2]+1 &\hbox{ if }t_2 \hbox{ is not an integer}.\end{cases}
\end{aligned}
\end{equation}
(Here $[t]$ stands for the integer part of $t$ and $t_1,t_2$ are as in \eqref{t_1}.) We have
\begin{equation*}
\begin{aligned}
&v^d_{\a}(x_0,y_0)=(1-\a)\left(\sum_{t=0}^{\t_1}\a^t+\sum_{t=\t_2}^{\infty}\a^t\right)=1-\a^{\t_1+1}+\a^{\t_2}
=1-\a^{t_1}\a^{p_1}+\a^{t_2}\a^{p_2}\\
\end{aligned}
\end{equation*}
where $p_1=\t_1+1-t_1$, $p_2=\t_2-t_2$. Observe that $0\le p_1,p_2< 1$.
Therefore
\begin{equation}\label{BB2}
v^d_{\a}(x_0,y_0)=1-\a^{t_1}+\a^{t_2}+\psi_{\a}
\end{equation}
where 
\begin{equation*}
\psi_{\a}(x_0,y_0)=\a^{t_1}(1-\a^{p_1})-\a^{t_2}(1-\a^{p_2}).
\end{equation*}
One can see that $\psi_{\a}\to 0$ as $\a\up 1$ uniformly with respect to $p_i\in [0,1)$ and $t_i\in [0,\infty)$, $i=1,2$. ($p_i$ and $t_i$ depend on $(x_0,y_0)$.)

By substitution $s=\a^{1/y_0}$, we obtain $1-\a^{t_1}+\a^{t_2}=1-s^{1-x_0}(1-s)$. Due to the calculations after formula \eqref{min_v}, for $\a\in [1/2,1)$ we have
\begin{equation}\label{X8}
\min_{y_0\in (0,1]}(1-\a^{t_1}+\a^{t_2})=1-{(1-x_0)^{1-x_0}\o (2-x_0)^{2-x_0}}.
\end{equation}
(Condition $\a\ge 1/2$ ensures that the point of minimum $s=(1-x_0)/(2-x_0)$ is in the interior of the interval $[0,\a]$ for all $x_0\in [0,1]$.)

Denote $\disp \nu_1(\a):=\sup_{(x_0,y_0)\in D}|\psi_{\a}(x_0,y_0)|$. Then, from \eqref{BB2} and \eqref{X8} we obtain
\begin{equation}\label{X5}
\left|\inf_{y_0\in (0,1]}v^d_{\a}(x_0,y_0)-
\left(1-{(1-x_0)^{1-x_0}\o (2-x_0)^{2-x_0}}\right)\right|\le\nu_1(\a)\,.
\end{equation}
Since $y_0=0$ is not a point of minimum of $v^d_{\a}(x_0,\cdot)$, $\disp\inf_{y_0\in (0,1]}v^d_{\a}(x_0,y_0)$ in \eqref{X5} can be replaced with $\disp\inf_{y_0\in [0,1]}v^d_{\a}(x_0,y_0)$, leading to  \eqref{VV5.7}.

\hf

\bigskip

It follows from \eqref{ZX1} and \eqref{VV5.7} that
\begin{equation}\label{VV5.6}
\begin{aligned}
&\inf_{y_0\in [0,1]}\g^d_T(0,y_0)= {1\o 2}+\eta_1(T),\\
&\inf_{y_0\in [0,1]}v^d_{\a}(0,y_0)={3\o 4}+\nu_1(\a)\,.
\end{aligned}
\end{equation}

\subsection{Smoothening}

Let $\ve>0$ be a sufficiently small parameter to be specified later and let $g^{\ve}:\,[0,\infty)\to \reals$ be a piecewise-linear  continuous function given by 
$$
g^{\ve}(x)=\begin{cases} 
1,&0\le x\le 1,\\ \disp 1-{x-1\o \ve}, &1<x<1+\ve,\\ 
0,&1+\ve\le x\le 2-\ve,\\ \disp {x-2+\ve\o \ve},&2-\ve<x<2,\\
1,& x\ge 2\,. \end{cases}
$$

On the trajectories of \eqref{VV4} denote 
\begin{equation*}
\begin{aligned}
&\g^{d,\ve}_T(x_0,y_0)={1\o T}\sum_{t=0}^{T-1} g^{\ve}(x(t)),\\
&v^{d,\ve}_{\a}(x_0,y_0)=(1-\a)\sum_{t=0}^{\infty} \a^{t}g^{\ve}(x(t)). 
\end{aligned}
\end{equation*}

{\bf Proposition 5.} {\em For $x_0\in [0,1)$ we have
\begin{equation}\label{ZX13}
\min_{y_0\in [0,1]}\g^{d,\ve}_T(x_0,y_0)= {1-x_0\o 2-x_0}+\eta_2(T)+\om_1(\ve),
\end{equation}
where $\eta_2(T)\to 0$ as $T\to \infty$ uniformly in $\ve$ and $x_0\in[0,1)$,  $\om_1(\ve)\to 0$ as $\ve\dn 0$ uniformly in $T$ and $x_0\in[0,1)$. 
}

{\bf Proof.} For $x_0\in [0,1)$ and $y_0\in (0,1]$, let $\t_1'\in N$ be the maximum time when the trajectory of \eqref{VV4} satisfies $x(\t_1')\le 1+\ve$ and $\t_2'\in N$ be the minimum time when this trajectory satisfies $x(\t_2')\ge 2-\ve$. As in Proposition 4,  let $\t_2\in N$ be the minimum time when this trajectory satisfies $x(\t_2)\ge 2$.
We have
\begin{equation}\label{NZ1}
\begin{aligned}
\g^{d,\ve}_T(x_0,y_0)&\le   {1\o T}\left(\sum_{t=0}^{\t'_1} 1+\sum_{t=\t_2'}^{\t_2-1}1+\max(T-\t_2,0)\right)\\
&\le{1\o T}\left(({1+\ve-x_0\o y_0}+1)+({\ve\o y_0}+1)+\max(T-\t_2,0)\right).
\end{aligned}
\end{equation}
Let $y_0^*=(2-x_0)/T$ to be the minimizer of $\g^{c}_T(x_0,y_0)$ with respect to $y_0\in [0,1]$, as in Proposition 1. Then
$\t_2=T$ since $x(\t_2)=2$ and, from \eqref{NZ1},
\begin{equation*}
\begin{aligned}
\inf_{y_0\in (0,1]}\g^{d,\ve}_T(x_0,y_0)\le\g^{d,\ve}_T(x_0,y^*_0)
={1-x_0\o 2-x_0}+{2\ve\o  2-x_0}+{2\o T}< {1-x_0\o 2-x_0}+{2\ve}+{2\o T} .
\end{aligned}
\end{equation*}
Combining this with the fact $y_0=0$ is not a point of minimum of $\g^{d,\ve}_T(x_0,\cdot)$, and that, due to \eqref{ZX1}, 
$$
\min_{y_0\in [0,1]}\g^{d,\ve}_T(x_0,y_0)\ge \inf_{y_0\in [0,1]}\g^{d}_T(x_0,y_0)={1-x_0\o 2-x_0}+\eta_1(T),
$$
we arrive at \eqref{ZX13} with $\eta_2(T)=\eta_1(T)+2/T$ and 
$\disp \om_1(\ve)={2\ve}$.

\bigskip

{\bf Proposition 6.} {\em For $x_0\in [0,1)$ and $\a\in [1/2,1)$ we have
\begin{equation}\label{ZX14}
\min_{y_0\in [0,1]}v^{d,\ve}_{\a}(x_0,y_0)=1-{(1-x_0)^{1-x_0}\o (2-x_0)^{2-x_0}}+\nu_2(\a)+\om_2(\ve)\,,
\end{equation}
where $\om_2(\ve)\to 0$ as $\ve\dn 0$ uniformly in $\a$ and $x_0\in[0,1)$, and $\nu_2(\a)\to 0$ as $\a\up 1$ uniformly in $\ve$ and $x_0\in[0,1)$.
}

{\bf Proof.}  For $x_0\in [0,1)$, $y_0\in (0,1]$, $\t_1'$ and $\t_2'$ as in the previous proposition, we have
\begin{equation}\label{X6}
\begin{aligned}
&v^{d,\ve}_{\a}(x_0,y_0)\le (1-\a)\left(\sum_{t=0}^{\t_1'}\a^t+\sum_{t=\t'_2}^{\infty}\a^t\right)
=1-\a^{\t'_1+1}+\a^{\t'_2}
\le 1-\a^{((1+\ve-x_0)/y_0)+1}+\a^{(2-\ve-x_0)/y_0}\\
&=1-\a^{(1+\ve-x_0)/y_0}+\a^{(2-\ve-x_0)/y_0}+(1-\a)\a^{(1+\ve-x_0)/y_0}.
\end{aligned}
\end{equation}
Notice that the last term aproaches 0 as $\a\up 1$ uniformly with respect to $x_0\in[0,1)$, $y_0\in (0,1]$, and $\ve>0$.
Denoting $s=\a^{1/y_0}$ and analyzing $1-s^{1+\ve-x_0}+s^{2-\ve-x_0}$ for a minimum on $(0,\a]$ by taking the derivative and setting it equal to zero, we obtain that the minimum is reached for $s=(1+\ve-x_0)/(2-\ve-x_0)$ and is equal to 
\begin{equation*}
\begin{aligned}
&1-\left({1+\ve-x_0\o 2-\ve-x_0}\right)^{1+\ve-x_0}+\left({1+\ve-x_0\o 2-\ve-x_0}\right)^{2-\ve-x_0}
=1-\left({1-x_0\o 2-x_0}\right)^{1-x_0}+\left({1-x_0\o 2-x_0}\right)^{2-x_0}+\om_2(\ve)\\
&=1-{(1-x_0)^{1-x_0}\o (2-x_0)^{2-x_0}}+\om_2(\ve),
\end{aligned}
\end{equation*}
where 
$$
\om_2(\ve):=-\left({1+\ve-x_0\o 2-\ve-x_0}\right)^{1+\ve-x_0}+\left({1+\ve-x_0\o 2-\ve-x_0}\right)^{2-\ve-x_0}-\left({1-x_0\o 2-x_0}\right)^{1-x_0}+\left({1-x_0\o 2-x_0}\right)^{2-x_0}
$$
approaches 0 as $\ve\dn 0$ uniformly with respect to $x_0\in [0,1)$, and does not depend on $\a$.
Thus, we obtain from \eqref{X6} that 
\begin{equation}\label{X9}
\inf_{y_0\in (0,1]}v^{d,\ve}_{\a}(x_0,y_0)\le 
1-{(1-x_0)^{1-x_0}\o (2-x_0)^{2-x_0}}+\om_2(\ve)+\nu_3(\a),
\end{equation}
where 
$$
\nu_3(\a):=\sup_{x_0,y_0} (1-\a)\a^{(1-x_0)/y_0}
$$
and supremum is taken over $(x_0,y_0)\in [0,1]\times (0,1]$.
Combining \eqref{X9} with the fact $y_0=0$ is not a point of minimum of $v^{d,\ve}_{\a}(x_0,\cdot)$, and that, due to \eqref{VV5.7}, 
$$
\min_{y_0\in [0,1]}v^{d,\ve}_{\a}(x_0,y_0)\ge \min_{y_0\in [0,1]}v^{d,\ve}_{\a}(x_0,y_0)=1-{(1-x_0)^{1-x_0}\o (2-x_0)^{2-x_0}}+\nu_1(\a),
$$
we arrive at \eqref{ZX14} with $\nu_2(\a):=\max\{\nu_1(\a),\nu_3(\a)\}.$

\hf

\subsection{Discrete Time Controlled System}

Now consider a discrete controlled system similar to \eqref{VV4} where, for convenience, the initial condition is given at $t=-1$. 
\begin{equation}\label{disc}
\begin{aligned}
&x(t+1)=x(t)+y(t), \; t\in\{ -1,0,1,\dots\},\\
&y(t+1)= u(t),\\
&(x(-1),y(-1))=(x_0,y_0),
\end{aligned}
\end{equation}
where $u(t)\in U(x(t),y(t))$ is given by 
$$
U(x,y)=\begin{cases} \disp \left[{xy\o a},1-{x(1-y)\o a}\right],\,&0\le x \le a\\  \{y\},\,&x> a.\end{cases}
$$
Here $a\in (0,1)$ is a sufficiently small parameter to be specified later. 
Note that $U(\cdot,\cdot)$ satisfies the following properties:

(a) $U(0,y)=[0,1]$ for all $y\in [0,1]$;

(b) $U(x,y)=\{y\}$ for $x\ge a$ and $y\in [0,1]$;

(c) $U(\cdot,\cdot)$ is Hausdorff continuous on $D$.


Property (a) implies that, if $x(t)=0$, control $u(t)$ can be selected so that $y(t+1)$ is anywhere in $[0,1]$. From property (b) it follows that, if $x(t)\ge a$, the system is uncontrollable and its dynamics coincides with the dynamics of \eqref{VV4}.

Along with \eqref{disc} consider a system controllable only at the initial state:
 \begin{equation}\label{disc1}
\begin{aligned}
&\tilde x(t+1)=\tilde x(t)+\tilde y(t), \; t\in \{-1,0,1,\dots\},\\
&\tilde y(t+1)= v(t),\\
&(x(-1),y(-1))=(0,0),
\end{aligned}
\end{equation}
where  $v(t)\in \tilde U(\tilde x(t),\tilde y(t))$ and 
$$
\tilde  U(x,y)=\begin{cases} \disp [0,1],\,&(x,y)=(0,0),\\  \{y\},\,&\hbox{otherwise.}\end{cases}
$$
If $v(\tau)\neq 0$ for some $\tau$, then \eqref{disc1} is uncontrollable for $t>\tau$, and its dynamics coincides with the dynamics of \eqref{VV4}.

\bigskip

{\bf Proposition 7.} {\em The $x$-component of any admissible process  in \eqref{disc} with $(x_0,y_0)=(0,0)$ can be approximated by the $x$-component an admissible process in \eqref{disc1}. More precisely, if $(x(\cdot),y(\cdot))$ is admissible \eqref{disc} with $(x_0,y_0)=(0,0)$, there exists a process $(\tilde x(\cdot),\tilde y(\cdot))$ admissible in \eqref{disc1}, such that
\begin{equation}\label{MM3}
|x(t)-\tilde x(t)|\le a \quad \hbox{for all }t\ge -1.
\end{equation}
}

This proposition is proved in the Appendix.

\hf

\bigskip

On the trajectories of \eqref{disc} denote 
\begin{equation}\label{VV10}
\begin{aligned}
&V_T^{a,\ve}(x_0,y_0):=\min_{u(\cdot)} {1\o T}\sum_{t=0}^{T-1} g^{\ve}(x(t)),\\
&h_{\a}^{a,\ve}(x_0,y_0):=\min_{u(\cdot)} (1-\a)\sum_{t=0}^{\infty} \a^{t}g^{\ve}(x(t)). 
\end{aligned}
\end{equation}

{\bf Proposition 8.} {\em Functions $V_T^{a,\ve}(\cdot,\cdot)$ and $h_{\a}^{a,\ve}(\cdot,\cdot)$ are continuous on $[0,\infty)\times [0,1]$.}

This proposition is proved in the Appendix.

\bigskip

The main result is the following theorem.

{\bf Theorem 1.} {\em For sufficiently small $a$ and $\ve$
\begin{equation}\label{main_ineq}
\lim_{T\to \infty} V_T^{a,\ve}(0,0) < \lim_{\a\up 1} h_{\a}^{a,\ve}(0,0).
\end{equation}
}

\bigskip

{\bf Proof.} On the trajectories $(\tilde x(\cdot),\tilde y(\cdot))$ of \eqref{disc1} denote
\begin{equation}\label{VV14}
\begin{aligned}
&\tilde V_T^{\ve}:=\min_{v(\cdot)} {1\o T}\sum_{t=0}^{T-1} g^{\ve}(\tilde x(t)),\\
&\tilde h_{\a}^{\ve}:=\min_{v(\cdot)} (1-\a)\sum_{t=0}^{\infty} \a^{t}g^{\ve}(\tilde x(t)). 
\end{aligned}
\end{equation}
The optimal control in \eqref{disc1} with former criterion in \eqref{VV14} is
$v(t)\equiv y'$, where $y'$ is the minimizer in \eqref{ZX13} with $x_0=0$; the optimal control in \eqref{disc1} with the latter criterion in \eqref{VV14} is
$v(t)\equiv y''$, where $y''$ is the minimizer in \eqref{ZX14} with $x_0=0$. Due to Propositions 5 and 6,
\begin{equation}\label{WW1}
\begin{aligned}
&\tilde V_T^{\ve}=\min_{y\in [0,1]} \g_T^{d,\ve}(0,y)={1\o 2}+\eta_2(T)+\om_1(\ve),\\
&\tilde h_{\a}^{\ve}=\min_{y\in [0,1]} v_{\a}^{d,\ve}(0,y)={3\o 4}+\nu_2(\a)+\om_2(\ve).
\end{aligned}
\end{equation}

Since any trajectory admissible in  \eqref{disc1} with $(x_0,y_0)=(0,0)$ is also admissible in \eqref{disc}, we have 
$$
V_T^{a,\ve}(0,0)\le \tilde V_T^{\ve}.
$$
Let $(x(\cdot),y(\cdot))$ be optimal in \eqref{disc} with $(x_0,y_0)=(0,0)$ with the first criterion in \eqref{VV10}.
Due to Proposition 7  there exists a process $(\tilde x(\cdot),\tilde y(\cdot))$ admissible in \eqref{disc1}, such that \eqref{MM3} holds. For this process
\begin{equation*}
\begin{aligned}
V_T^{a,\ve}(0,0)&= {1\o T}\sum_{t=0}^{T-1} g^{\ve}(x(t))
 \ge {1\o T}\sum_{t=0}^{T-1} g^{\ve}(\tilde x(t))-{1\o T}\sum_{t=0}^{T-1} |g^{\ve}(\tilde x(t))-g^{\ve}(x(t))|.
\end{aligned}
\end{equation*}
Since $g^{\ve}(\cdot)$ is equicontinuous on $[0,\infty)$, we have $|g^{\ve}(x')-g^{\ve}(x'')|\to 0$ as $|x'-x''|\to 0$ uniformly with respect to $x',x''$. Therefore, inequality above implies that
$$
V_T^{a,\ve}(0,0)\ge \tilde V_T^{\ve}-\mu_1^{\ve}(a),
$$
where $\disp\mu_1^{\ve}(a):={1\o T}\sum_{t=0}^{T-1} |g^{\ve}(\tilde x(t))-g^{\ve}(x(t))|$ is a quantity that approaches zero as $a\dn 0$ uniformly with respect to $T$ due to \eqref{MM3}. 

Since any trajectory admissible in  \eqref{disc1} with $(x_0,y_0)=(0,0)$ is also admissible in \eqref{disc}, we have 
$$
h_{\a}^{a,\ve}(0,0)\le \tilde h_{\a}^{\ve}.
$$
Let $(x(\cdot),y(\cdot))$ be optimal in \eqref{disc} with $(x_0,y_0)=(0,0)$ with the second criterion in \eqref{VV10}. 
Similarly to the lines above, there exists a process $(\tilde x(\cdot),\tilde y(\cdot))$ admissible in \eqref{disc1}, such that \eqref{MM3} holds. For this process
$$
h_{\a}^{a,\ve}(0,0)=(1-\a)\sum_{t=0}^{\infty} \a^{t}g^{\ve}(x(t))\ge(1-\a)\sum_{t=0}^{\infty} \a^{t}g^{\ve}(\tilde x(t))-\mu_2^{\ve}(a)\ge
 \tilde h_{\a}^{\ve}-\mu_2^{\ve}(a),
$$
where $\mu_2^{\ve}(a)$ is a quantity that approaches zero as $a\dn 0$ uniformly with respect to $\a$. 
Thus
\begin{equation*}
\begin{aligned}
&\tilde V_T^{\ve}-\mu_1^{\ve}(a)\le V_T^{a,\ve}(0,0)\le  \tilde V_T^{\ve},\\
&\tilde h_{\a}^{\ve}-\mu_2^{\ve}(a)\le h_{\a}^{a,\ve}(0,0)\le \tilde h_{\a}^{\ve}\,.
\end{aligned}
\end{equation*}
Taking into account \eqref{WW1} we conclude from these inequalities above that
\begin{equation*}
\begin{aligned}
&|V_T^{a,\ve}(0,0)-{1\o 2}|\le \mu_1^{\ve}(a)+\eta_2(T)+\om_1(\ve)\\
&|h_{\a}^{a,\ve}(0,0)-{3\o 4}|\le \mu_2^{\ve}(a)+\nu_2(\a)+\om_2(\ve).
\end{aligned}
\end{equation*}
Therefore, $V_T^{a,\ve}(0,0)$ can be made arbitrarily close to 1/2 by taking $T$ sufficiently large, $\ve$ and $a$ sufficiently small; $h_{\a}^{a,\ve}(0,0)$ can be made arbitrarily close to 3/4 by taking $\a$ sufficiently close to 1, $\ve$ and $a$ sufficiently small. This implies the validity of Theorem 1.

\hf

\bigskip

In the constructed example \eqref{disc} the system evolves on the unbounded set $[0,\infty)\times [0,1]$. To be consistent with the assumptions of \cite{BGS} this set should be compact. This compactification can be achieved by changing the first equation in the dynamics of \eqref{disc} to
$x(t+1)=x(t)+y(t)+q(x(t))$, where $q(x)$ is a continuous decreasing function such that $q(x)=0$ for $0\le x\le 2$ and $q(x)=-x$ for $x\ge 3$. Then the system only evolves on a compact invariant set, and this modification doesn't effect the values of limits in \eqref{main_ineq}. 

Also notice that, consistently with \eqref{A22}, we have equality
$$
\lim_{T\to \infty}\min_{(x_0,y_0)}V_T^{a,\ve}(x_0,y_0)=\lim_{\a\up 1}\min_{(x_0,y_0)}h_{\a}^{a,\ve}(x_0,y_0),
$$
since both quantities are equal to zero for $x_0 \in (1+\ve,2-\ve)$ and $y_0=0$.

\section{Appendix}

{\bf Proof of Proposition 7.}

Let $(x(t),y(t)),\,t\ge -1$ be an admissible process in \eqref{disc} with $(x_0,y_0)=(0,0)$. If $x(t)\le a$ for all $t$, then \eqref{MM3} holds with $\tilde x\equiv 0$. Assume that  
$(x(\cdot),y(\cdot))$ is such that $x(\t)>a$ for some $\t$. Let $\bar t$ be the maximum moment of time such that $x(\bar t)\le a$. Since $x(0)=x(-1)=0$, it follows that $\bar t\ge 1$.

From \eqref{disc} and due to the facts that $y(t)$ is constant (see property (b) of $U$) for $t\ge \bar t+1$ and $(x_0,y_0)=(0,0)$, we have
\begin{equation}\label{MM1}
x(t)=\begin{cases}\disp \sum_{\t=0}^{t-1}y(\t),&t\le \bar t\\ \disp\sum_{\t=0}^{\bar t-1}y(\t)+(t-\bar t)y(\bar t),&t\ge \bar t+1\,.\end{cases}
\end{equation}
Let $v(\cdot)$ in \eqref{disc1} be
\begin{equation*}
v(t)=\begin{cases}\disp 0, &t\le \bar t-2\\  y(\bar t),&t\ge \bar t-1\,.\end{cases}
\end{equation*}
Then for the corresponding $(\tilde x(\cdot),\tilde y(\cdot))$ we have 
$$
\tilde y(t)=\begin{cases} 0,&t\le \bar t-1\\ y(\bar t),&t\ge \bar t\,\end{cases}
$$
and
\begin{equation}\label{MM2}
\tilde x(t)=\sum_{\t=0}^{t-1}y(\t)=\begin{cases} 0,&t\le \bar t\\ (t-\bar t)y(\bar t),&t\ge \bar t+1\,.\end{cases}
\end{equation}
Taking into account that $\disp x(\bar t)=\sum_{\t=0}^{\bar t-1}y(\t)\le a$ by the definition of $\bar t$, we conclude from \eqref{MM1} and \eqref{MM2} that
\begin{equation*}
0\le x(t)-\tilde x(t)\le a \quad \hbox{for all }t,
\end{equation*}
which is the required approximation property. 

\hf

{\bf Proof of Proposition 7.}
Assume that the data in \eqref{A1} is such that  the set $Y$ is invariant, that is, for any $y\in Y$ and $u\in U(y)$ we have $f(y,u)\in Y$. Also assume that $f$ and $g$ are uniformly continuous and bounded on $Y$, the mapping $U(\cdot)$ is compact-valued and uniformly Hausdorff continuous on $Y$. (The uniformity assumptions hold automatically when $Y$ is a compact set.) Let us show continuity of $V_T(\cdot)$ and $h_{\a}(\cdot)$ in \eqref{A112-1} and \eqref{A112-2}. This will imply that Proposition 7 is true, since the data of \eqref{disc} satisfies the assumptions above.

{\em Continuity of $V_T$.} 
Let $(y(\cdot),u(\cdot))$ be an optimal process in problem \eqref{A112-1}. 

From continuity of $f$ and Hausdorff continuity of $U(\cdot)$ it follows that for any $\ve>0$ there exists $\d>0$ such that for any $y_0'$ such that $|y_0'-y_0|<\d$ there exists $u'(0)\in U(y_0')$ such that $|f(y_0',u'(0))-f(y_0,u(0))|<\ve$. This implies that 
for $y'(1):= f(y_0',u'(0))$ we have
$$
|y'(1)-y(1)|= |f(y_0',u'(0))-f(y_0,u(0))|< \ve.
$$
Continuing this process we can construct an admissible trajectory $y'(t),\,t=0,\dots,T$  of \eqref{A1} such that $|y'(t)-y(t)|<M\ve$ for any $t$ with some constant $M$ that depends only on $T$. Therefore, for any $\gamma>0$ and any $y_0'$ sufficiently close to $y_0$ we have
$$
V_T(y_0')-V_T(y_0)\le{1\o T}\sum_{t=0}^{T-1} g(y'(t))-{1\o T}\sum_{t=0}^{T-1}g(y(t))\le {1\o T}\sum_{t=0}^{T-1} |g(y'(t))-g(y(t))|\le \gamma.
$$
Hence, $|V_T(y_0')-V_T(y_0)| \le \gamma$, that is, $V_T(\cdot)$ is continuous.

{\em Continuity of $h_{\a}(\cdot)$.} 
Let $(y(\cdot),u(\cdot))$ be an optimal process in \eqref{A112-2} (an optimal process exists due to Proposition 1 in \cite{GPS-2017}. For any admissible process  $(y'(t),u'(t)),\,t=0,\dots,T$  we have
\begin{equation*}
\begin{aligned}
&V_T(y_0')-V_T(y_0)\le (1-\a)\sum_{t=0}^{k} \a^t|g(y'(t))-g(y(t))|+(1-\a)\sum_{t=k+1}^{\infty} \a^t|g(y'(t))-g(y(t))|.
\end{aligned}
\end{equation*}
Take $\gamma>0$. Take $k$ sufficiently large so that $(1-\a)\sum_{t=k+1}^{\infty} \a^t|g(y'(t))-g(y(t))|\le \gamma/2$. Then for any $y_0'$ sufficiently close to $y_0$ there exists an admissible trajectory $y'(t),\,t=0,\dots,T$ so that $y'(0)=y_0'$ and the first sum does not exceed $\gamma/2$. In this case  
 $|V_T(y_0')-V_T(y_0)| \le \gamma$.

\hf


\bigskip





\begin{thebibliography}{99}

 \bibitem{AN} E. J. Anderson and P. Nash, Linear Programming in Infinite-Dimensional Spaces, Wiley, Chichester, 1987.

 \bibitem{Arisawa-3} M. Arisawa and  P.-L. Lions, On Ergodic Stochastic Control, {\it Commun. in Partial Differential Equations}, 23:11 (1998), pp.  2187--2217.



\bibitem{Bardi} M. Bardi and I. Capuzzo-Dolcetta, Optimal Control and Viscosity Solutions of Hamilton-Jacobi-Bellman Equations, Birkhauser, Boston, 1997.



\bibitem{BGS}
V. Borkar, V. Gaitsgory and I. Shvartsman,
LP Formulations of  Discrete Time Long-Run Average Optimal Control Problems: The Non-Ergodic Case, {\em SIAM Journal on Control and Optimization}, 57:3 (2019), pp. 1783--1817.


\bibitem{BQR-2015} R. Buckdahn, M. Quincampoix and J. Renault,  On Representation Formulas for Long Run Averaging Optimal control Problem, {\it Journal of Differential Equations}, 259:11 (2015), pp. 5554--5581.



\bibitem{GQ} V. Gaitsgory and M. Quincampoix,
Linear Programming Approach to Deterministic
Infinite Horizon Optimal Control Problems with Discounting,  {\it SIAM J. of Control and
Optimization},  48:4 (2009), pp. 2480-2512.

\bibitem {GPS-2017} V. Gaitsgory, A. Parkinson and I. Shvartsman, Linear Programming Formulations of Deterministic Infinite Horizon Optimal Control Problems in Discrete Time, {\it Discrete and Continuous Dynamical Systems Series B},  22:10 (2017), pp. 3821--3338.






\bibitem{Gonzalez98} O. Hernandez-Lerma and J. Gonzales-Hernandez, Infinite Linear Programming and Multichain
Markov Control Processes in Uncountable Spaces, {\em SIAM J. of Control and Optimization}, 36:1
(1998), pp. 313-335.

\bibitem{GruneSIAM98}
 L. Gr\"une,
 Asymptotic Controllability and Exponential Stabilization of Nonlinear Control Systems at Singular Points,
 {\it SIAM J. of Control and Optimization}, 36:5 (1998), pp. 1495--1503.



\bibitem{GruneJDE98}
L. Gr\"une,
 On the Relation Between Discounted and Average Optimal Value Functions.
{\it J. Diff. Equations}, 148 (1998), pp. 65--69.


\bibitem{HardyLit14}  G.H. Hardy, J.E. Littlewood, 
Tauberian theorems concerning power series and 
Dirichlet series whose coefficients are positive.
Proc. London Math. Soc. 13, pp. 174--191 (1914) 

\bibitem{Hardy49}  G.H. Hardy, Divergent series. Clarendon Press, Oxford, 1949. 

\bibitem{Khlopin}  D. Khlopin, Tauberian Theorem for Value Functions, {\it Dynamical Games and Applications}, 8 (2018), pp. 401--422.


\bibitem{Sorin92}
 E. Lehrer and S.  Sorin,
  A Uniform Tauberian Theorem in Dynamic Programming,
{\it Mathematics of Operations Research}, 17:2,  (1992), pp. 303--307.




\bibitem{OV-2012}  M. Oliu-Barton and G.   Vigeral,   A Uniform Tauberian Theorem in Optimal Control, In  P. Cardaliaguet and R.  Grossman (Eds.){\it Annals of International Society of Dynamic Games},  12, pp. 199-215, Birkhauser/Springer, New York, 2013.


\bibitem{QR-2011} Quincampoix, M., Renault, J.: On the existence of a limit value in some non expansive optimal control problems, {\em SIAM J. Control Optim.} 49 (2011), pp. 2118--2132. 



\end{thebibliography}
\end{document}